\documentclass[english]{elsarticle}
\usepackage{ae,aecompl}
\usepackage[LGR,T1]{fontenc}
\usepackage[latin9]{inputenc}
\usepackage{array}
\usepackage{booktabs}
\usepackage{textcomp}
\usepackage{multirow}
\usepackage{amsmath}
\usepackage{graphicx}
\usepackage{esint}
\usepackage{subscript}

\usepackage{geometry}
\makeatletter

\journal{Composite Structures}

\DeclareRobustCommand{\greektext}{%
  \fontencoding{LGR}\selectfont\def\encodingdefault{LGR}}
\DeclareRobustCommand{\textgreek}[1]{\leavevmode{\greektext #1}}
\DeclareFontEncoding{LGR}{}{}
\DeclareTextSymbol{\~}{LGR}{126}
\providecommand{\tabularnewline}{\\}

\usepackage{lscape}

\@ifundefined{showcaptionsetup}{}{%
 \PassOptionsToPackage{caption=false}{subfig}}
\usepackage{subfig}
\makeatother

\usepackage{babel}
\begin{document}

\title{\textbf{Smoothed Finite Element and Genetic Algorithm based optimization
for Shape Adaptive Composite Marine Propellers}}

\author{Manudha T Herath$^{1*}$}

\author{Sundararajan Natarajan$^{2}$}

\author{B Gangadhara Prusty$^{1}$}

\author{Nigel St John$^{3}$}

\address{$^{1}$School of Mechanical and Manufacturing Engineering,University
of New South Wales, Sydney, Australia}

\address{$^{2}$School of Civil and Environmental Engineering, University
of New South Wales, Sydney, Australia}

\address{$^{3}$Maritime Division, Defence Science and Technology Organization,
Fishermans Bend, Victoria, Australia}

\address{{*}Corresponding author Email: m.herath@unsw.edu.au, Tel.: +61293855054}
\begin{abstract}
An optimization scheme using the Cell-based Smoothed Finite Element
Method (CS-FEM) combined with a Genetic Algorithm (GA) framework is
proposed in this paper to design shape adaptive laminated composite
marine propellers. The proposed scheme utilise the bend-twist coupling
characteristics of the composites to achieve the required performance.
An iterative procedure to evaluate the unloaded shape of the propeller
blade is proposed, confirming the manufacturing requirements at the
initial stage. The optimization algorithm and codes developed in this
work were implemented under a variety of parameter settings and compared
against the requirement to achieve an ideally passive pitch varying
propeller. Recommendations for the required thickness of the propeller
blade to achieve optimum bend-twist coupling performance without resulting
large rake deformations are also presented.\end{abstract}
\maketitle

\section{Introduction}

Marine propellers are traditionally manufactured from Nickel Aluminium
Bronze (NAB) alloys or Manganese Nickel Aluminium Bronze (MAB) alloys.
However, recently the use of engineered materials, more specifically
laminated composites, to develop marine propellers has received considerable
attention equally among researchers and industry. This is driven by
the increasing demand for efficiency, high strength-to-weight and
high stiffness-to-weight ratio of materials. Some of the advantages
exhibited by composite over alloys are light weight, reduced corrosion,
reduced noise generation, lack of magnetic signature and shape adaptability.
Shape-adaptability is an interesting phenomenon from a mechanical
and optimization perspective. 

Shape-adaptability refers to the capability of composites to deform,
without the involvement of external mechanisms based on the flow conditions
and rotational speed in order to achieve a higher efficiency, compared
to rigid alloy propellers, throughout its operating domain. Composite
propellers can potentially be custom tailored to enhance the performance
through shape-adaptability. This can be achieved by using the intrinsic
extension-shear, bend-twist and bend-extension coupling effects of
anisotropic composites. Ideal shape change at various flow conditions
is a result of composite layup optimizations such that the propeller
has an optimum bend-twist coupling performance.

Bend-twist coupling refers to the special characteristic of anisotropic
materials where out of plane bending moments can cause twisting strains.
With correct layup arrangements this effect can be optimized for a
certain application using layered composites. Various researchers
in the past have used flexibility and bend-twist coupling characteristics
of composites to design marine propellers that have the capability
of self-varying pitch (shape adaptable) based on out of plane bending
moments caused by the incoming flow \citep{LeeLin_2004,young2007,Liu2009,Mulcahy_Offshore,Motley_etal_2011,MulcahyRINA_2011}. 

The approach taken by Lin and Lee \citep{LeeLin_2004,Linetal_2009,Linetal2004}
was to minimize the change of torque coefficient of the propeller
when moving from the design advance ratio to one other off-design
advance ratio. The reason behind this strategy was maintaining the
torque, thrust and efficiency the same as the design value when moving
away from the design point. However, only one off-design point was
considered. The optimization process used previously by Liu, et al.
\citep{Liu2009}, Motley, et al. \citep{Motley_etal_2011}, Pluci\'{n}ski,
et al. \citep{Plucinski2007} attempted to ensure that the ply configuration
was chosen such that the blade can achieve the maximum possible pitch
variation when moving from unloaded to loaded state. Essentially,
the optimization technique attempted to make the blade more flexible
while maintaining strain and shape limitations. 

In this proposed Finite Element (FE) approach, the propeller blade
assuming it is plate idealized to the mid-plane of the overall blade
shape. Various structural theories proposed for evaluating the characteristics
of composite laminates under different loading situations were reviewed
by Noor and Burton (1989) \citep{noor1989}, Mallikarjuna and Kant
(1993) \citep{mallikarjunakant1993}, Kant and Swaminathan (2000)
\citep{kantswaminathan2000} and recently by Khanda et. al \citep{khandannoroozi2012}.
A set of methods have emerged to address the shear locking in the
FEM. By incorporating the strain smoothing technique into the finite
element method (FEM), Liu et al. \citep{liudai2007} have formulated
a series of smoothed finite element methods (SFEM), named as cell-based
SFEM (CS-FEM) \citep{nguyenbordas2008,bordasnatarajan2010}, node-based
SFEM \citep{liunguyen2009}, edge-based SFEM \citep{liunguyen2009a},
face-based SFEM \citep{thoiliu2009} and \textgreek{a}-FEM \citep{liunguyen2008}.
And recently, edge based imbricate finite element method (EI-FEM)
was proposed in \citep{Cazes2012} that shares common features with
the ES-FEM. As the SFEM can be recast within a Hellinger-Reissner
variational principle, suitable choices of the assumed strain/gradient
space provides stable solutions. Depending on the number and geometry
of the sub-cells used, a spectrum of methods exhibiting a spectrum
of properties is obtained. Further details can be found in other literature
\citep{nguyenbordas2008} and references therein. Nguyen-Xuan et al.
\citep{nguyenrabczuk2008} employed CS-FEM for Mindlin-Reissner plates.
The curvature at each point is obtained by a non-local approximation
via a smoothing function. From the numerical studies presented, it
was concluded that the CS-FEM technique is robust, computationally
inexpensive, free of locking and importantly insensitive to mesh distortions.
The SFEM was extended to various problems such as shells \citep{nguyenrabczuk2008},
heat transfer \citep{wuliu2010}, fracture mechanics \citep{nguyenliu2013}
and structural acoustics \citep{hecheng2011} among others. In \citep{bordasnatarajan2011},
CS-FEM has been combined with the extended FEM to address problems
involving discontinuities.

A framework to design laminated composite marine propellers with enhanced
performance by utilizing the bend-twist coupling characteristics is
proposed in this paper. The framework consists of the Cell-based Smoothed
Finite Element Method (CS-FEM) combined with a Genetic Algorithm (GA)
to optimize the layup configuration of laminated composites. The key
requirement for the optimization technique proposed here is to achieve
an efficiency curve for the composite propeller that is tangential
to all the efficiency curves in the vicinity of the design (cruise)
advance ratio of the vessel. In contrast to the approaches taken by
previous researchers \citep{Leeetal2005,Plucinski2007,Liu2009}, the
proposed method attempts to achieve accurate pitch angles derived
from propeller efficiency curves based on many off-design points.
It also gives the freedom to specify weightages to off-design points
based on probabilities the blade is likely to operate at each off-design
point. The approach was presented by means of a simple plate optimisation
study by Herath and Prusty \citep{Herath_ACAM7}. It was further enhanced
to an Iso-Geometric (NURBS) FEM based optimisation technique by Herath
et al. \citep{Herath2013}.

\section{Theoretical Development}

The shape-adaptive technique presented in this paper predominantly
relies on bend-twist coupling characteristics of laminated composites
to change the pitch of the blade based on bending caused by fluid
loadings at different flow speeds. Bend-twist coupling characteristics
can be demonstrated using the standard stiffness matrix system for
composite materials (eq. \ref{eq:ABBD}). Here, $\begin{array}{cc}
\left[A\right], & \left[B\right]\end{array}$ and $\left[D\right]$ matrices have the usual laminate stiffness
definition. 

\begin{equation}
\begin{array}{c}
\left[\begin{array}{c}
\mathbf{N}\\
\mathbf{M}
\end{array}\right]=\left[\begin{array}{cc}
\left[A\right] & \left[B\right]\\
\left[B\right] & \left[D\right]
\end{array}\right]\left[\mathbf{\begin{array}{c}
\boldsymbol{\boldsymbol{\mathbf{\mathbf{\epsilon}}}}\\
\boldsymbol{\kappa}
\end{array}}\right]\;\textup{\textup{where;}}\\
\mathbf{N}=\mathrm{\mathrm{\left[\begin{array}{ccc}
N_{x} & N_{y} & N_{xy}\end{array}\right]}^{T},}\:\mathbf{M}=\mathrm{\left[\begin{array}{ccc}
M_{x} & M_{y} & M_{xy}\end{array}\right]^{T}}\\
\mathbf{\boldsymbol{\epsilon}=\mathrm{\left[\begin{array}{ccc}
\epsilon_{x} & \epsilon_{y} & \epsilon_{xy}\end{array}\right]^{T}}\mathrm{,}}\:\mathbf{\boldsymbol{\kappa}=\mathrm{\left[\begin{array}{ccc}
\kappa_{x} & \kappa_{y} & \kappa_{xy}\end{array}\right]^{T}}}
\end{array}\label{eq:ABBD}
\end{equation}

Bend-twist coupling characteristics are dominated through coupling
terms ($D_{xs}$ and $D_{ys}$) in the $\left[D\right]$ matrix (eq.
\ref{eq:Dmat}), where $D_{ij}=\frac{1}{3}\sum_{k=1}^{n}Q_{ij}^{k}\left(\theta\right)\left(z_{k}^{3}-z_{k-1}^{3}\right)$
with $Q_{ij}\left(\theta\right)$ representing the in-plane stiffness
of a composite layer in the $xy$ directions. With $D_{xs},\: D_{ys}\neq0$
bending moments ($M_{xx}$ and $M_{yy}$) can cause twisting strains
$\left(\frac{\partial^{2}w}{\partial x\partial y}\right)$. The purpose
of an optimization scheme is to obtain the optimum fiber angles to
achieve the required response (pressure vs required pitch). The matrix
system gives the relationship for a simple laminate element. For a
plate-like structure the stiffness coefficients have to be used in
combination with plate theories such as Kirchhoff-Love (thin plate)
and Mindlin-Reissner (moderately thick plate) with appropriate boundary
conditions. Thus, it is essential to have an FEM technique that can
accurately determine the response of a complex blade shape. 

\begin{equation}
\begin{array}{c}
\left[\begin{array}{c}
M_{xx}\\
M_{yy}\\
M_{xy}
\end{array}\right]\end{array}=\left[\begin{array}{ccc}
D_{xx} & D_{xy} & D_{xs}\\
D_{xy} & D_{yy} & D_{ys}\\
D_{xs} & D_{ys} & D_{ss}
\end{array}\right]\left[\begin{array}{c}
\kappa_{xx}\\
\kappa_{yy}\\
\kappa_{xy}
\end{array}\right]=\left[\begin{array}{ccc}
D_{xx} & D_{xy} & D_{xs}\\
D_{xy} & D_{yy} & D_{ys}\\
D_{xs} & D_{ys} & D_{ss}
\end{array}\right]\left[\begin{array}{c}
-\frac{\partial^{2}w}{\partial x^{2}}\\
-\frac{\partial^{2}w}{\partial y^{2}}\\
-2\frac{\partial^{2}w}{\partial x\partial y}
\end{array}\right]\label{eq:Dmat}
\end{equation}

\subsection{Genetic Algorithm based optimization}

The key idea behind the proposed optimization scheme for a marine
propeller is to construct a \textquotedblleft{}difference-scheme\textquotedblright{}
relative to the operating point (cruise advance ratio) in terms of
pressure and twist. The optimum alloy propeller geometry must be chosen
for the application before it is further developed as a composite
propeller. The process can be summarized as: 
\begin{enumerate}
\item Evaluate pressure maps on the propeller blade surface for various
speeds including and around the operating/cruise speed. 
\item Construct pressure difference functions with respect to the operating
condition for every chosen point around the operating point. 
\item Assess the pitch changes required relative to the operating point
for the chosen points to maintain an optimum efficiency. Pitch differences
can be assessed using standard propeller efficiency curves for a propeller
series, which the alloy propeller is based upon.
\item The objective function of optimization will attempt to minimize the
total difference (corresponding to the respective pressure difference)
between the optimum pitch that is required and the pitch that was
obtained by the chosen ply configuration (Eq. \ref{eq:ObjFunc}).
Weightages $\left(w_{i}\right)$ can be assigned to each off design
point based on the likelihood of the propeller operating at each off-design
point.
\begin{equation}
\begin{array}{c}
\begin{array}{c}
min\\
\mathbf{\boldsymbol{\theta}}
\end{array}f\left(\mathbf{\boldsymbol{\theta}}\right)=\frac{\sum_{i=1}^{n}w_{i}\left|\Delta\phi_{Optimum}^{i}\left(\Delta P\right)-\Delta\phi_{GA}^{i}\left(\mathbf{\boldsymbol{\theta}},\Delta P\right)\right|}{\sum_{i=1}^{n}w_{i}}\\
\textup{s.t. failure criteria and strain limitations}
\end{array}\label{eq:ObjFunc}
\end{equation}

\end{enumerate}
The optimization technique must be capable of handling non-linear
objective functions, non-linear constraints and both discrete and
continuous variables. Thus, the Genetic Algorithm (GA) was chosen
as it can satisfy all these requirements. The GA has been used by
several authors \citep{Kameyamaetal2007,Leeetal2005,Linetal2004,Plucinski2007,SoremekunGurdal2001}
in composite ply optimization tasks proving its attractiveness and
credibility. 

\begin{figure}
\begin{centering}
\includegraphics[scale=0.5]{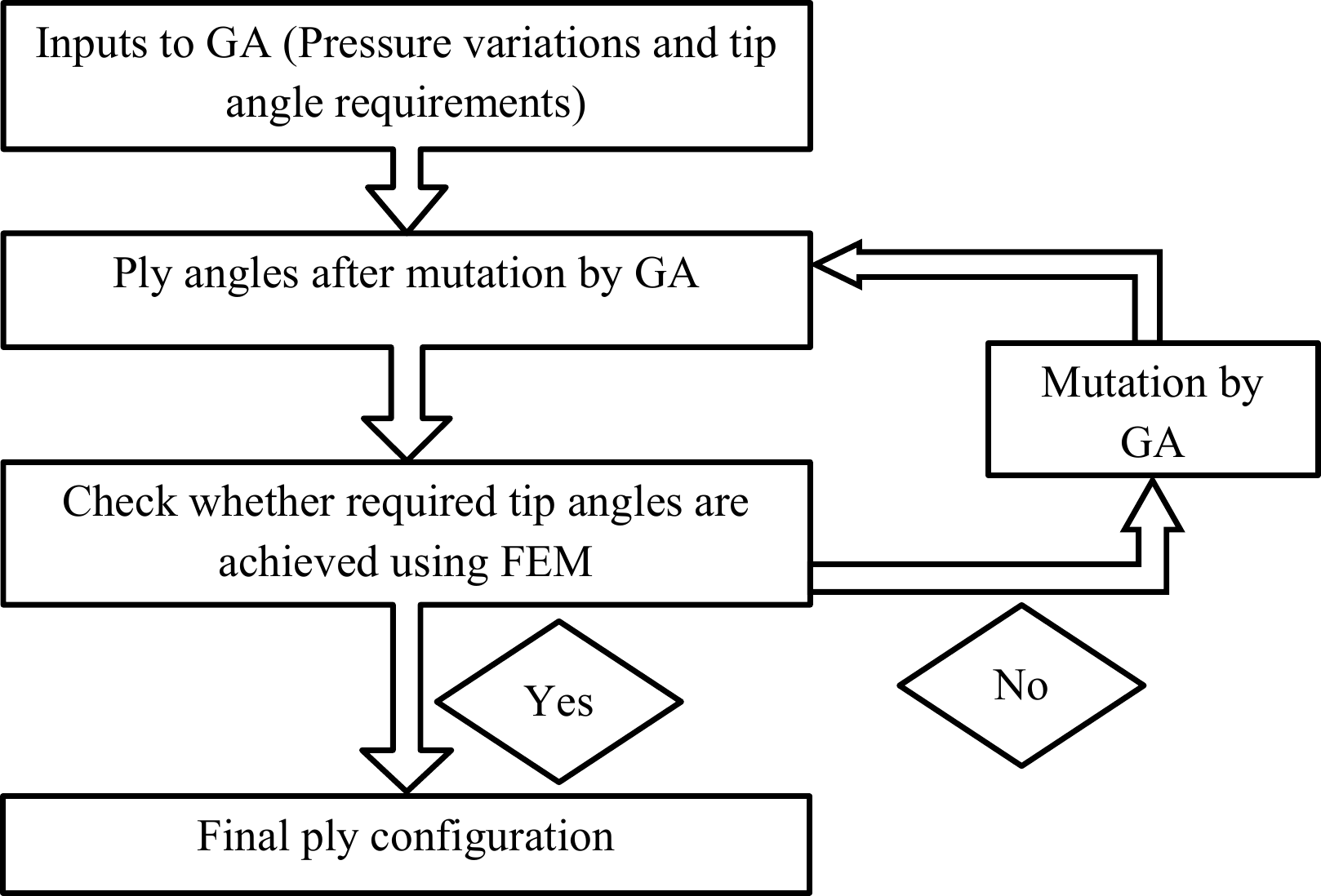}
\par\end{centering}

\caption{\label{fig:FEM-coupled-with}FEM coupled with GA flow}
\end{figure}

The process of GA involves applying mutations to the ply angle configuration
and evaluating whether the blade can achieve the required angle at
the tip. This gives rise to the requirement of having an accurate
means of calculating deflections and rotations of the blade structure
for an applied loading. Thus, an in-house FE code based on the first
order shear theory using Cell Based Smoothed FEM was developed for
propeller blade shapes and coupled with the GA. Figure \ref{fig:FEM-coupled-with}
shows a summary of the optimization process coupled with FEM. Both
the GA and the FE solver were coded in the commercial numerical processing
software Matlab\texttrademark{}. Although it is possible to couple
the GA with an existing commercial FEM solver as attempted by several
authors in similar research \citep{LeeLin_2004,Leeetal2005,Motley_etal_2011,MulcahyRINA_2011,Plucinski2007},
a coupled fully in-house solver is seen as a future proof approach.
This is due to the inherent freedom the user has within such a solver
and capability of improvement and further streamlining in the future. 

Furthermore, a composite propeller optimized for pitch variation cannot
be manufactured at its optimum shape required at the cruise speed.
The blade has to be manufactured with a pre-deformation such that
it reaches the optimum shape at its cruise condition. Thus, the second
stage of the design process involves an iterative process to evaluate
the unloaded shape. The proposed methodology is iterative as explained
by an example in Section \ref{sub:Unloaded-shape-calculation}. A
similar methodology was also used by Mulcahy, et al. \citep{MulcahyRINA_2011}
and Pluci\'{n}ski, et al. \citep{Plucinski2007}.

\subsection{Cell based smoothed finite element method with discrete shear gap
technique}

\label{csdsg3}In this study, the propeller blade is approximated
by a hypothetical plate at the mid-plane of the blade. Three-noded
triangular element with five degrees of freedom (dofs) $\boldsymbol{\delta}=\left\{ u,v,w,\theta_{x},\theta_{y}\right\} $
is employed to discretise the plate domain. The displacement is approximated
by 
\begin{equation}
\mathbf{u}^{h}=\sum_{I}N_{I}\boldsymbol{\delta}_{I}
\end{equation}

where $\boldsymbol{\delta}_{I}$ are the nodal dofs and $N_{I}$ are
the standard finite element shape functions given by 
\begin{equation}
N=\left[1-\xi-\eta,\;\;\eta,\;\;\xi\right]
\end{equation}

In the CS-DSG3~\cite{thoivan2012}, each triangular element is divided into three sub-triangles.
The displacement vector at the center node is assumed to be the simple
average of the three displacement vectors of the three field nodes.
In each sub-triangle, the stabilized DSG3 is used to compute the strains
and also to avoid the transverse shear locking. Then the strain smoothing
technique on the whole triangular element is used to smooth the strains
on the three sub-triangles.

\begin{center}
\begin{figure}[htpb]
\begin{centering}
\includegraphics[scale=0.5]{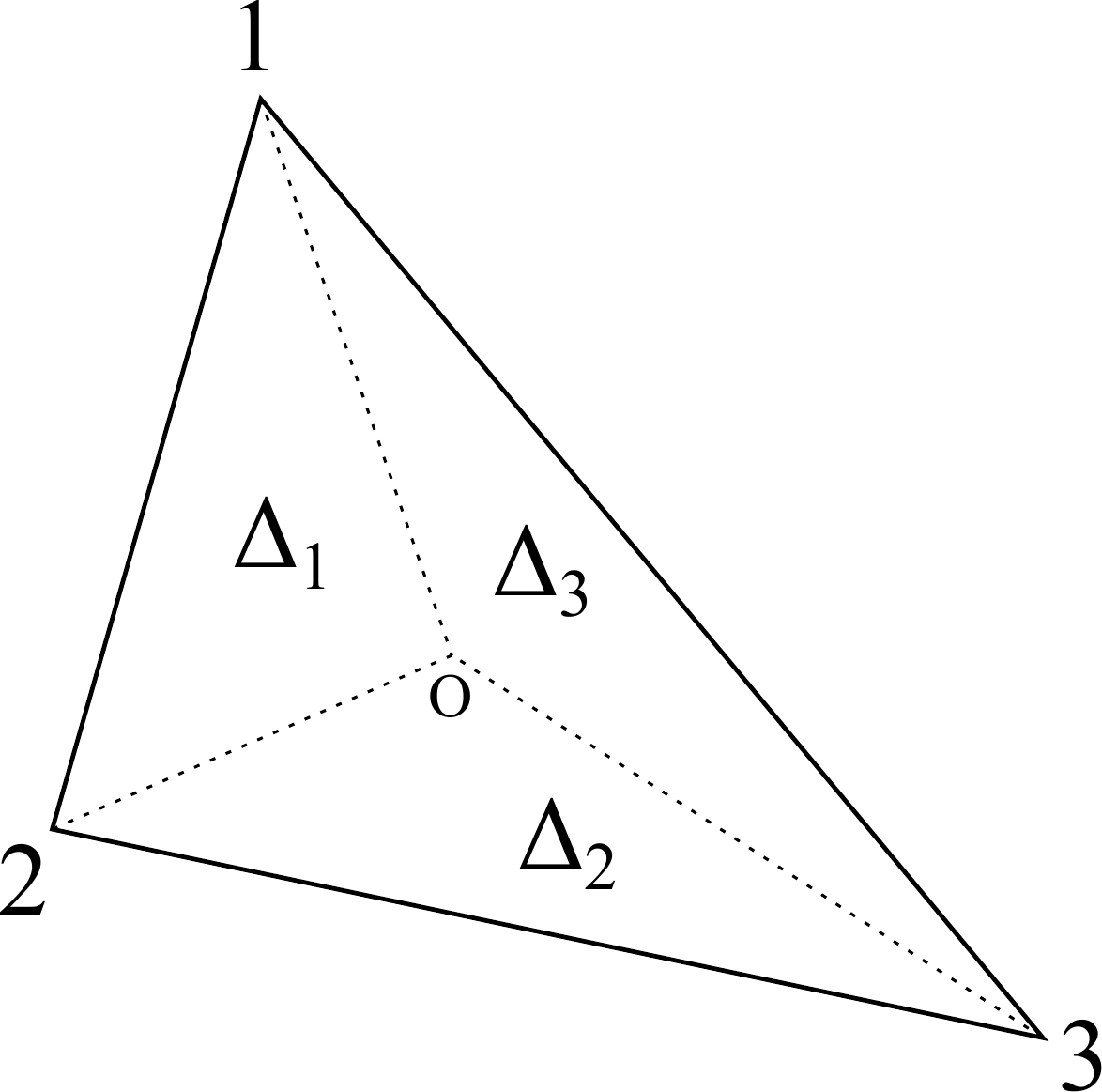}
\par\end{centering}

\caption{A triangular element is divided into three sub-triangles. $\Delta_{1},\Delta_{2}$
and $\Delta_{3}$ are the sub-triangles created by connecting the
central point $O$ with three field nodes.}

\label{fig:triEle} 
\end{figure}

\par\end{center}

Consider a typical triangular element $\Omega_{e}$ as shown in Figure
\ref{fig:triEle}. This is first divided into three sub-triangles
$\Delta_{1},\Delta_{2}$ and $\Delta_{3}$ such that $\Omega_{e}=\bigcup\limits _{i=1}^{3}\Delta_{i}$.
The coordinates of the center point $\mathbf{x}_{o}=(x_{o},y_{o})$
is given by: 
\begin{equation}
(x_{o},y_{o})=\frac{1}{3}(x_{I},y_{I})
\end{equation}
The displacement vector of the center point is assumed to be a simple
average of the nodal displacements as 
\begin{equation}
\boldsymbol{\delta}_{eO}=\frac{1}{3}\boldsymbol{\delta}_{eI}\label{eqn:centerdefl}
\end{equation}

The constant membrane strains, the bending strains and the shear strains
for sub-triangle $\Delta_{1}$ is given by:

\begin{equation}
\begin{array}{c}
\boldsymbol{\epsilon}_{p}=\left[\begin{array}{ccc}
\mathbf{p}_{1}^{\Delta_{1}} & \mathbf{p}_{2}^{\Delta_{1}} & \mathbf{p}_{3}^{\Delta_{1}}\end{array}\right]\left\{ \begin{array}{c}
\boldsymbol{\delta}_{eO}\\
\boldsymbol{\delta}_{e1}\\
\boldsymbol{\delta}_{e2}
\end{array}\right\} \\
\boldsymbol{\epsilon}_{b}=\left[\begin{array}{ccc}
\mathbf{b}_{1}^{\Delta_{1}} & \mathbf{b}_{2}^{\Delta_{1}} & \mathbf{b}_{3}^{\Delta_{1}}\end{array}\right]\left\{ \begin{array}{c}
\boldsymbol{\delta}_{eO}\\
\boldsymbol{\delta}_{e1}\\
\boldsymbol{\delta}_{e2}
\end{array}\right\} \\
\boldsymbol{\epsilon}_{s}=\left[\begin{array}{ccc}
\mathbf{s}_{1}^{\Delta_{1}} & \mathbf{s}_{2}^{\Delta_{1}} & \mathbf{s}_{3}^{\Delta_{1}}\end{array}\right]\left\{ \begin{array}{c}
\boldsymbol{\delta}_{eO}\\
\boldsymbol{\delta}_{e1}\\
\boldsymbol{\delta}_{e2}
\end{array}\right\} 
\end{array}\label{eq:constrains}
\end{equation}

Upon substituting the expression for $\boldsymbol{\delta}_{eO}$ in
Eqs. \ref{eq:constrains}, we obtain:

\begin{equation}
\begin{array}{c}
\boldsymbol{\epsilon}_{p}^{\Delta_{1}}=\left[\begin{array}{ccc}
\frac{1}{3}\mathbf{p}_{1}^{\Delta_{1}}+\mathbf{p}_{2}^{\Delta_{1}} & \frac{1}{3}\mathbf{p}_{1}^{\Delta_{1}}+\mathbf{p}_{3}^{\Delta_{1}} & \frac{1}{3}\mathbf{p}_{1}^{\Delta_{1}}\end{array}\right]\left\{ \begin{array}{c}
\boldsymbol{\delta}_{eO}\\
\boldsymbol{\delta}_{e1}\\
\boldsymbol{\delta}_{e2}
\end{array}\right\} =\mathbf{B}_{p}^{\Delta_{1}}\boldsymbol{\delta}_{e}\\
\boldsymbol{\epsilon}_{b}^{\Delta_{1}}=\left[\begin{array}{ccc}
\frac{1}{3}\mathbf{b}_{1}^{\Delta_{1}}+\mathbf{b}_{2}^{\Delta_{1}} & \frac{1}{3}\mathbf{b}_{1}^{\Delta_{1}}+\mathbf{b}_{3}^{\Delta_{1}} & \frac{1}{3}\mathbf{b}_{1}^{\Delta_{1}}\end{array}\right]\left\{ \begin{array}{c}
\boldsymbol{\delta}_{eO}\\
\boldsymbol{\delta}_{e1}\\
\boldsymbol{\delta}_{e2}
\end{array}\right\} =\mathbf{B}_{b}^{\Delta_{1}}\boldsymbol{\delta}_{e}\\
\boldsymbol{\epsilon}_{s}^{\Delta_{1}}=\left[\begin{array}{ccc}
\frac{1}{3}\mathbf{s}_{1}^{\Delta_{1}}+\mathbf{s}_{2}^{\Delta_{1}} & \frac{1}{3}\mathbf{s}_{1}^{\Delta_{1}}+\mathbf{s}_{3}^{\Delta_{1}} & \frac{1}{3}\mathbf{s}_{1}^{\Delta_{1}}\end{array}\right]\left\{ \begin{array}{c}
\boldsymbol{\delta}_{eO}\\
\boldsymbol{\delta}_{e1}\\
\boldsymbol{\delta}_{e2}
\end{array}\right\} =\mathbf{B}_{s}^{\Delta_{1}}\boldsymbol{\delta}_{e}
\end{array}
\end{equation}

where $\mathbf{p}_{i},\,(i=1,2,3),\:\mathbf{b}_{i},\,(i=1,2,3)$ and
$\mathbf{s}_{i},\,(i=1,2,3)$ are given by:

\begin{equation}
\begin{array}{c}
\mathbf{B}_{p}=\frac{1}{2A_{e}}\left[\underbrace{\begin{array}{c}
b-c\\
0\\
d-a
\end{array}}_{\mathbf{p}_{1}}\underbrace{\begin{array}{c}
0\\
d-a\\
-d
\end{array}}_{\mathbf{p}_{2}}\underbrace{\begin{array}{c}
0\\
0\\
a
\end{array}}_{\mathbf{p}_{3}}\begin{array}{rrrrrrrrrrrr}
0 & 0 & c & 0 & 0 & 0 & 0 & -b & 0 & 0 & 0 & 0\\
0 & 0 & 0 & -d & 0 & 0 & 0 & a & 0 & 0 & 0 & 0\\
0 & 0 & 0 & 0 & 0 & 0 & 0 & 0 & 0 & 0 & 0 & 0
\end{array}\right]\\
\mathbf{B}_{b}=\frac{1}{2A_{e}}\left[\underbrace{\begin{array}{c}
0\\
0\\
0
\end{array}}_{\mathbf{b}_{1}}\underbrace{\begin{array}{c}
0\\
0\\
0
\end{array}}_{\mathbf{b}_{2}}\underbrace{\begin{array}{c}
0\\
0\\
0
\end{array}}_{\mathbf{b}_{3}}\begin{array}{rrrrrrrrrrrr}
b-c & 0 & 0 & 0 & 0 & c & 0 & 0 & 0 & 0 & -b & 0\\
0 & d-a & 0 & 0 & 0 & 0 & -d & 0 & 0 & 0 & 0 & a\\
0 & 0 & 0 & 0 & 0 & 0 & 0 & 0 & 0 & 0 & 0 & 0
\end{array}\right]\\
\mathbf{B}_{s}=\frac{1}{2A_{e}}\left[\underbrace{\begin{array}{c}
0\\
0
\end{array}}_{\mathbf{s}_{1}}\underbrace{\begin{array}{c}
0\\
0
\end{array}}_{\mathbf{s}_{2}}\underbrace{\begin{array}{c}
b-c\\
0
\end{array}}_{\mathbf{s}_{3}}\begin{array}{rrrrrrrrrrrr}
A_{e} & 0 & 0 & 0 & c & \frac{ac}{2} & \frac{bc}{2} & 0 & 0 & -b & \frac{-bd}{2} & \frac{-bc}{2}\\
0 & 0 & 0 & 0 & 0 & 0 & 0 & 0 & 0 & 0 & 0 & 0
\end{array}\right]
\end{array}
\end{equation}

where $a=x_{2}-x_{1};b=y_{2}-y_{1};c=y_{3}-y_{1}$ and $d=x_{3}-x_{1}$
(see Figure \ref{fig:dsg3}), $A_{e}$ is the area of the triangular
element and $\mathbf{B}_{s}$ is altered shear strains \citep{bletzingerbischoff2000}.
The strain-displacement matrix for the other two triangles can be
obtained by cyclic permutation.

\begin{center}
\begin{figure}[htpb]
\begin{centering}
\includegraphics[scale=0.5]{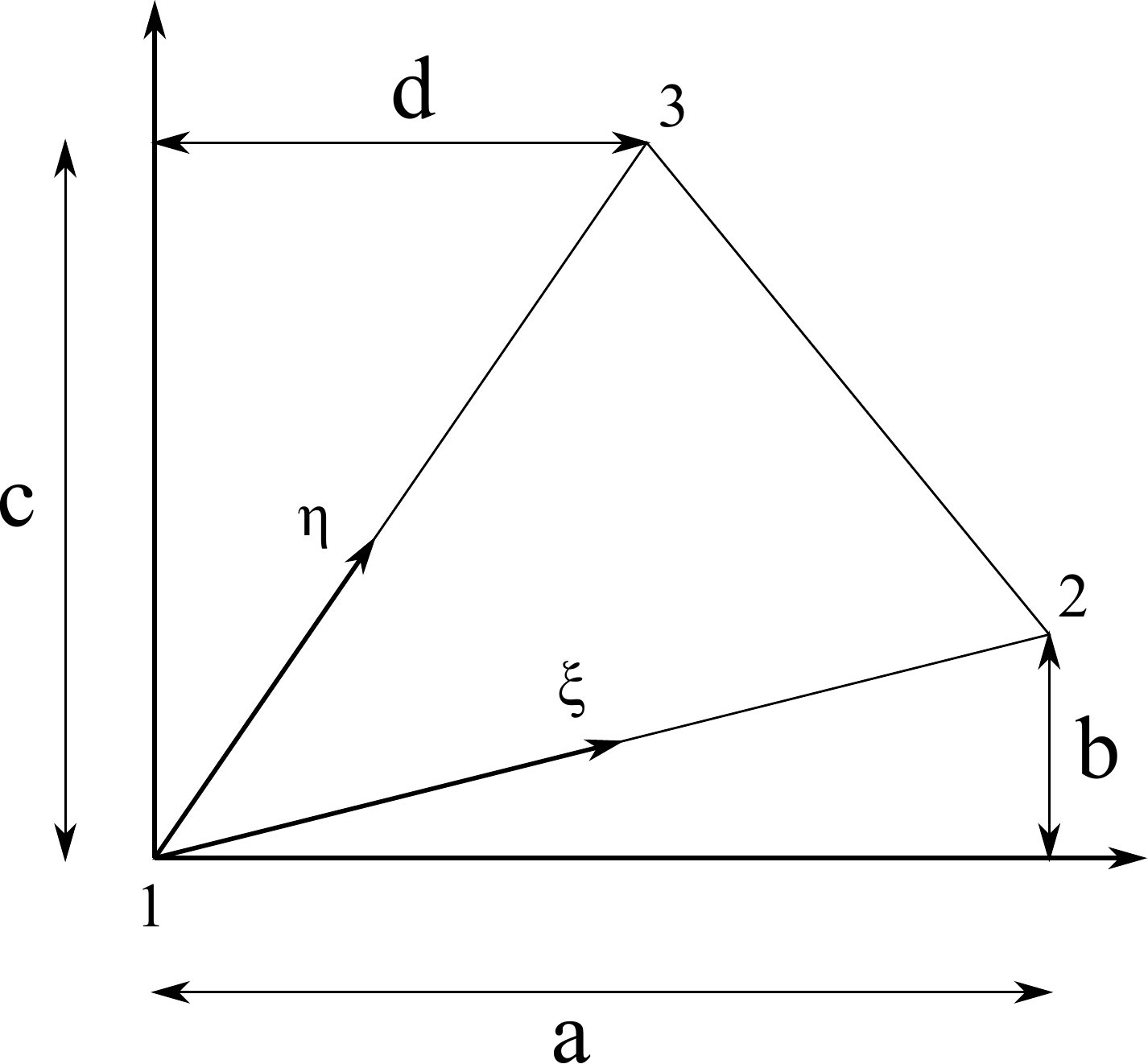}
\par\end{centering}

\caption{Three-noded triangular element and local coordinates in discrete shear
gap method.}

\label{fig:dsg3} 
\end{figure}

\par\end{center}

Now applying the cell-based strain smoothing \citep{bordasnatarajan2010},
the constant membrane strains, the bending strains and the shear strains
are respectively employed to create a smoothed membrane strain $\boldsymbol{\bar{\epsilon}}_{p}$,
smoothed bending strain $\boldsymbol{\bar{\epsilon}}_{b}$ and smoothed
shear strain $\boldsymbol{\bar{\epsilon}}_{s}$ on the triangular
element $\Omega_{e}$ as:

\begin{equation}
\begin{array}{c}
\boldsymbol{\bar{\epsilon}}_{p}=\int\limits _{\Omega_{e}}\boldsymbol{\epsilon}_{p}\Phi_{e}\left(\mathbf{x}\right)d\Omega=\sum\limits _{i=1}^{3}\boldsymbol{\epsilon}_{p}^{\Delta_{i}}\int\limits _{\Delta_{i}}\Phi_{e}(\mathbf{x})d\Omega\\
\boldsymbol{\bar{\epsilon}}_{b}=\int\limits _{\Omega_{e}}\boldsymbol{\epsilon}_{b}\Phi_{e}\left(\mathbf{x}\right)d\Omega=\sum\limits _{i=1}^{3}\boldsymbol{\epsilon}_{p}^{\Delta_{i}}\int\limits _{\Delta_{i}}\Phi_{e}(\mathbf{x})d\Omega\\
\boldsymbol{\bar{\epsilon}}_{s}=\int\limits _{\Omega_{e}}\boldsymbol{\epsilon}_{s}\Phi_{e}\left(\mathbf{x}\right)d\Omega=\sum\limits _{i=1}^{3}\boldsymbol{\epsilon}_{s}^{\Delta_{i}}\int\limits _{\Delta_{i}}\Phi_{e}(\mathbf{x})d\Omega
\end{array}
\end{equation}

where $\Phi_{e}\left(\mathbf{x}\right)$ is a given smoothing function
that satisfies. In this study, following constant smoothing function
is used: 
\begin{equation}
\Phi\left(\mathbf{x}\right)=\left\{ \begin{array}{cc}
1/A_{c} & \mathbf{x}\in\Omega_{c}\\
0 & \mathbf{x}\notin\Omega_{c}
\end{array}\right.
\end{equation}

where $A_{c}$ is the area of the triangular element, the smoothed
membrane strain, the smoothed bending strain and the smoothed shear
strain is then given by 
\begin{equation}
\left\{ \boldsymbol{\bar{\epsilon}}_{p},\boldsymbol{\bar{\epsilon}}_{b},\boldsymbol{\bar{\epsilon}}_{s}\right\} =\frac{\sum\limits _{i=1}^{3}A_{\Delta_{i}}\left\{ \boldsymbol{\epsilon}_{p}^{\Delta_{i}},\boldsymbol{\epsilon}_{b}^{\Delta_{i}},\boldsymbol{\epsilon}_{s}^{\Delta_{i}}\right\} }{A_{e}}
\end{equation}

The smoothed elemental stiffness matrix is given by 
\begin{align}
\mathbf{K} & =\int\limits _{\Omega_{e}}\overline{\mathbf{B}}_{p}\mathbf{A}\overline{\mathbf{B}}_{p}^{{\rm T}}+\mathbf{\overline{B}}_{p}\mathbf{B}\mathbf{\overline{B}}_{b}^{{\rm T}}+\mathbf{\overline{B}}_{b}\mathbf{B}\mathbf{\overline{B}}_{p}^{{\rm T}}+\mathbf{\overline{B}}_{b}\mathbf{D}\mathbf{\overline{B}}_{b}^{{\rm T}}+\mathbf{\overline{B}}_{s}\mathbf{E}\mathbf{\overline{B}}_{s}^{{\rm T}}d\Omega\nonumber \\
 & =\left(\overline{\mathbf{B}}_{p}\mathbf{A}\overline{\mathbf{B}}_{p}^{{\rm T}}+\mathbf{\overline{B}}_{p}\mathbf{B}\mathbf{\overline{B}}_{b}^{{\rm T}}+\mathbf{\overline{B}}_{b}\mathbf{B}\mathbf{\overline{B}}_{p}^{{\rm T}}+\mathbf{\overline{B}}_{b}\mathbf{D}\mathbf{\overline{B}}_{b}^{{\rm T}}+\mathbf{\overline{B}}_{s}\mathbf{E}\mathbf{\overline{B}}_{s}^{{\rm T}}\right)A_{e}
\end{align}
where $\overline{\mathbf{B}}_{p},\overline{\mathbf{B}}_{b}$ and $\overline{\mathbf{B}}_{s}$
are the smoothed strain-displacement matrix.

\section{Numerical Results}

Numerical results were obtained using GA coupled FEM solver for Wageningen-B
series propellers. Convergence and stability of the FEM solver was
validated using a simple plate example. The solver was then successfully
used to optimize standard propeller blade shapes. To ensure that the
solver and the optimization technique is applicable for various blade
shapes, three blades with different Expanded Area Ratios (EAR) were
considered. The B5-45 blade was further analyzed to optimize for a
weighted off-design point case, integer ply optimization and investigate
the effect of layer thickness. The unloaded shape was also calculated
for the B5-45 blade.

\subsection{CS-FEM mesh convergence and stability}

A mesh convergence study was conducted to ensure that the cell-based
smoothed finite element technique is stable and provides accurate
results with the increase in the number of degrees of freedom. The
mesh was refined by increasing the node number of the test structure
(h-refinement). A simple rectangular plate with dimensions: 0.4 m
(L) x 0.2 m (W) x 3 mm (t) was considered for the convergence study.
It was assumed that the plate was made out of unidirectional CFRP
(Table \ref{tab:Material-Properties}) and has 24 plies all having
a fiber orientation of 40\textdegree{} counter clockwise from x-axis
towards y-axis. The plate was assumed to be clamped at the left edge
and a uniform pressure loading (normal to the surface) of 100 Pa (upwards)
was applied on the top surface. Details of the meshes that were validated
and their results are given in Table 1. As an independent verification,
maximum deflection obtained using Q8 elements (using the commercial
software ANSYS\texttrademark{}, 8-noded shell 281) is also presented.
Convergence results showed that CS-FEM was highly accurate with good
stability and convergence. Thus, it can be used for complex shapes
in further applications.

\begin{table}
\begin{centering}
\begin{tabular}{ccc}
\toprule 
 & Node Array & Max. Deflection (mm)\tabularnewline
\midrule
\midrule 
Mesh 1 & 5\texttimes{}5 & 4.296 \tabularnewline
\midrule
\midrule 
Mesh 2 & 10\texttimes{}10 & 5.704\tabularnewline
\midrule
\midrule 
Mesh 3 & 20\texttimes{}20 & 6.087\tabularnewline
\midrule
\midrule 
Mesh 4 & 40\texttimes{}40 & 6.165\tabularnewline
\midrule
\midrule 
Mesh 5 & 80\texttimes{}80 & 6.194\tabularnewline
\midrule
\midrule 
\multicolumn{2}{c}{ANSYS\texttrademark{} Q8 (9841 Nodes)} & 6.212\tabularnewline
\bottomrule
\end{tabular}
\par\end{centering}

\caption{\label{tab:Mesh-convergence}Mesh convergence of CS-FEM}
\end{table}

\begin{figure}
\begin{centering}
\includegraphics[scale=0.5]{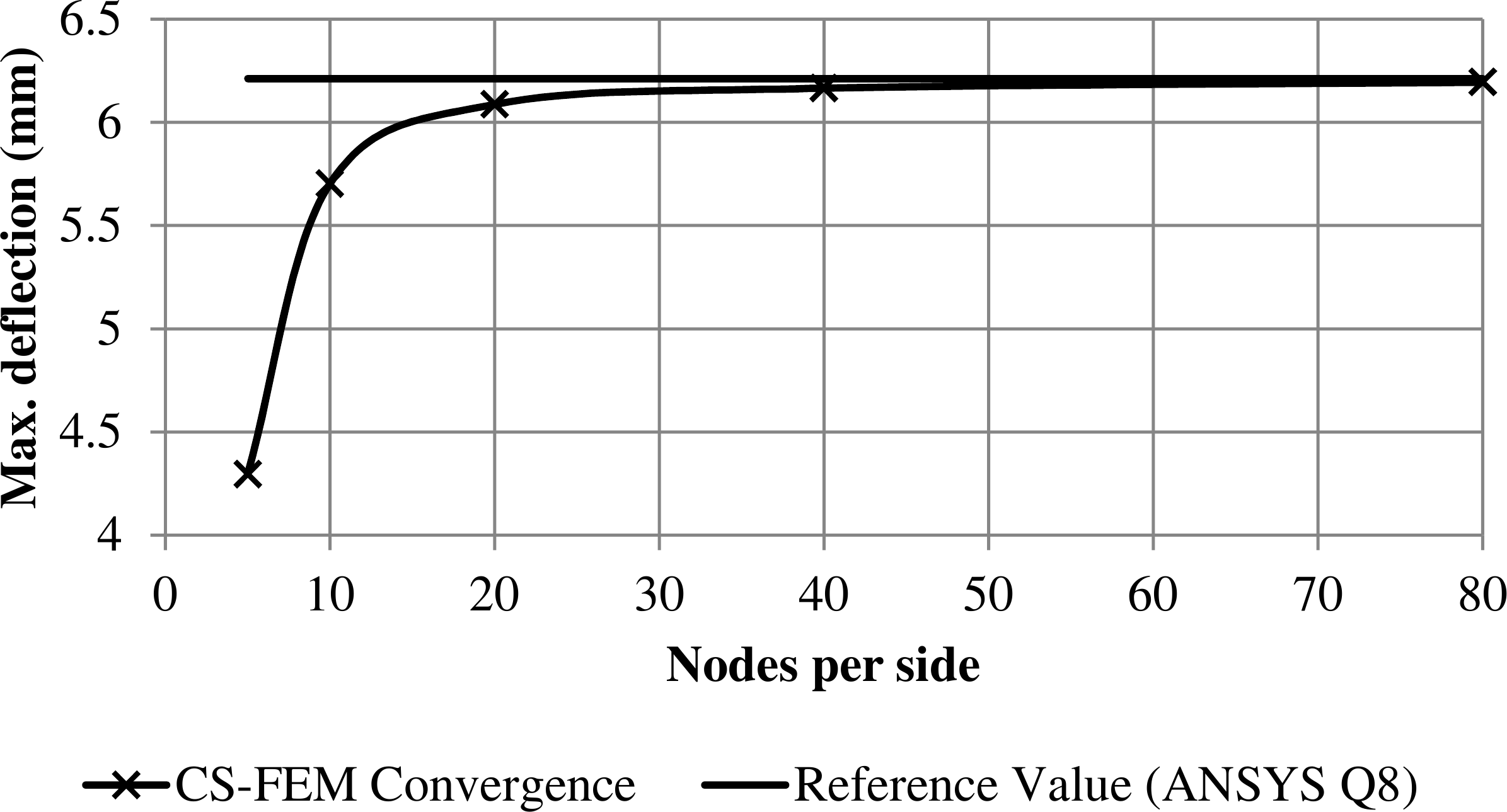}
\par\end{centering}

\caption{\label{fig:Mesh-convergence-curve}CS-FEM mesh convergence curve}

\end{figure}

\subsection{Preliminary layup optimization of three blade shapes}

The developed GA coupled FE solver can be used to optimize the Wageningen-B
Series propeller blades against various flow conditions at various
speeds. The proposed idea is to extract pressure distribution maps
at various speeds above and below the operating (cruise) advance ratio
and use the pressure maps as the basis of optimization. Method by
which pressure maps are evaluated is at the discretion of the hydrodynamist.
These pressure maps have to be used in the optimization scheme's finite
element routine as pressure differences $\left(\Delta P_{ij}\right)$
at each Gauss point calculated with respect to the pressure map at
cruise condition. However, this paper focuses on presenting the optimization
methodology rather than actual values; thus, pressure values are directly
used as the basis for optimization. Pressure variations used in this
paper were chosen as uniform arbitrary distributions with no obvious
relationship in order to maintain the generality. They were chosen
sensibly based on the pressure distribution at cruise condition obtained
using CFD (ANSYS CFX) analyzes in preliminary studies (Table \ref{tab:Pressure-variations}). 

\begin{figure}
\begin{centering}
\subfloat[]{\centering{}\includegraphics[scale=0.5]{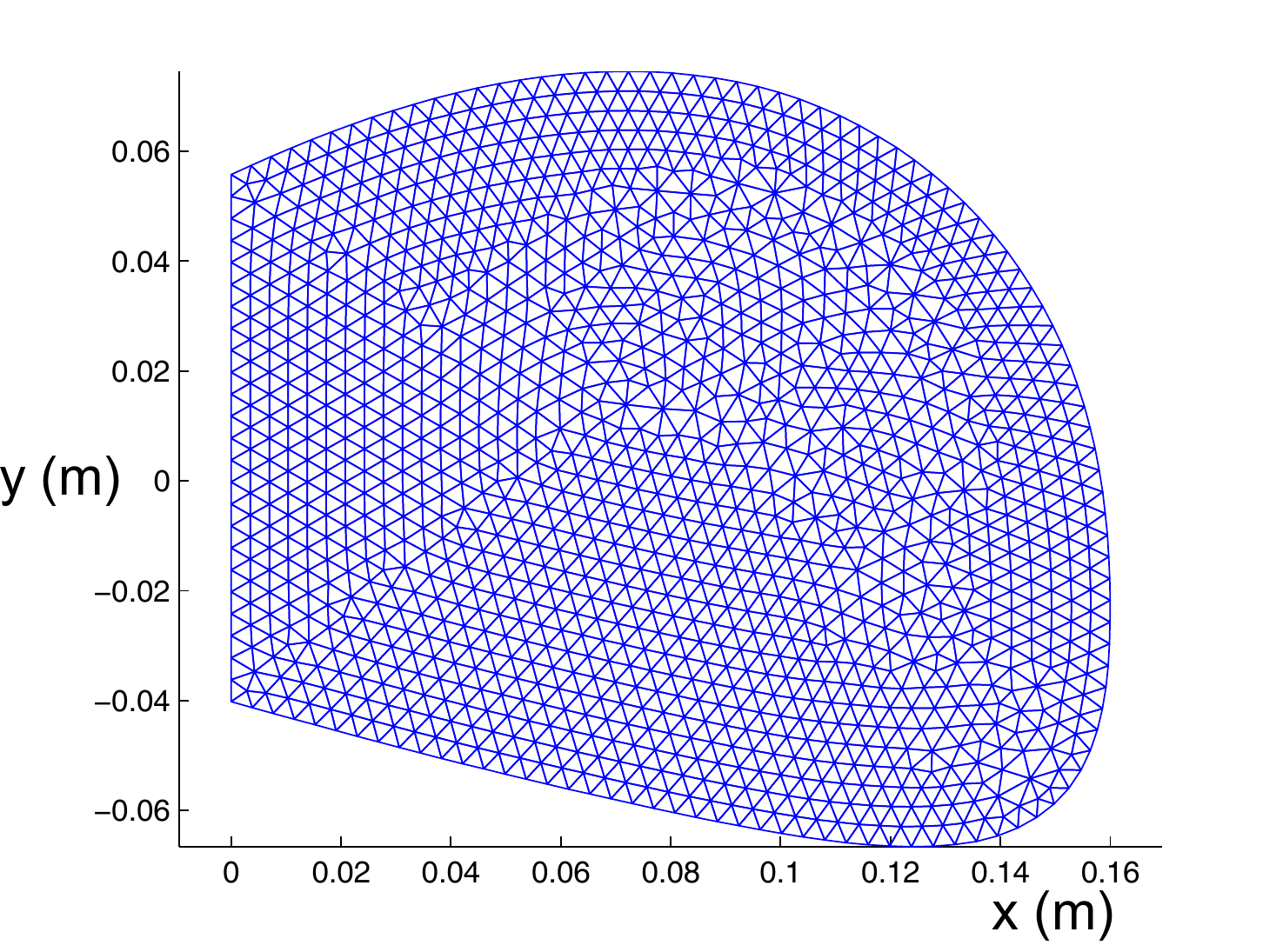}}\subfloat[]{\centering{}\includegraphics[scale=0.5]{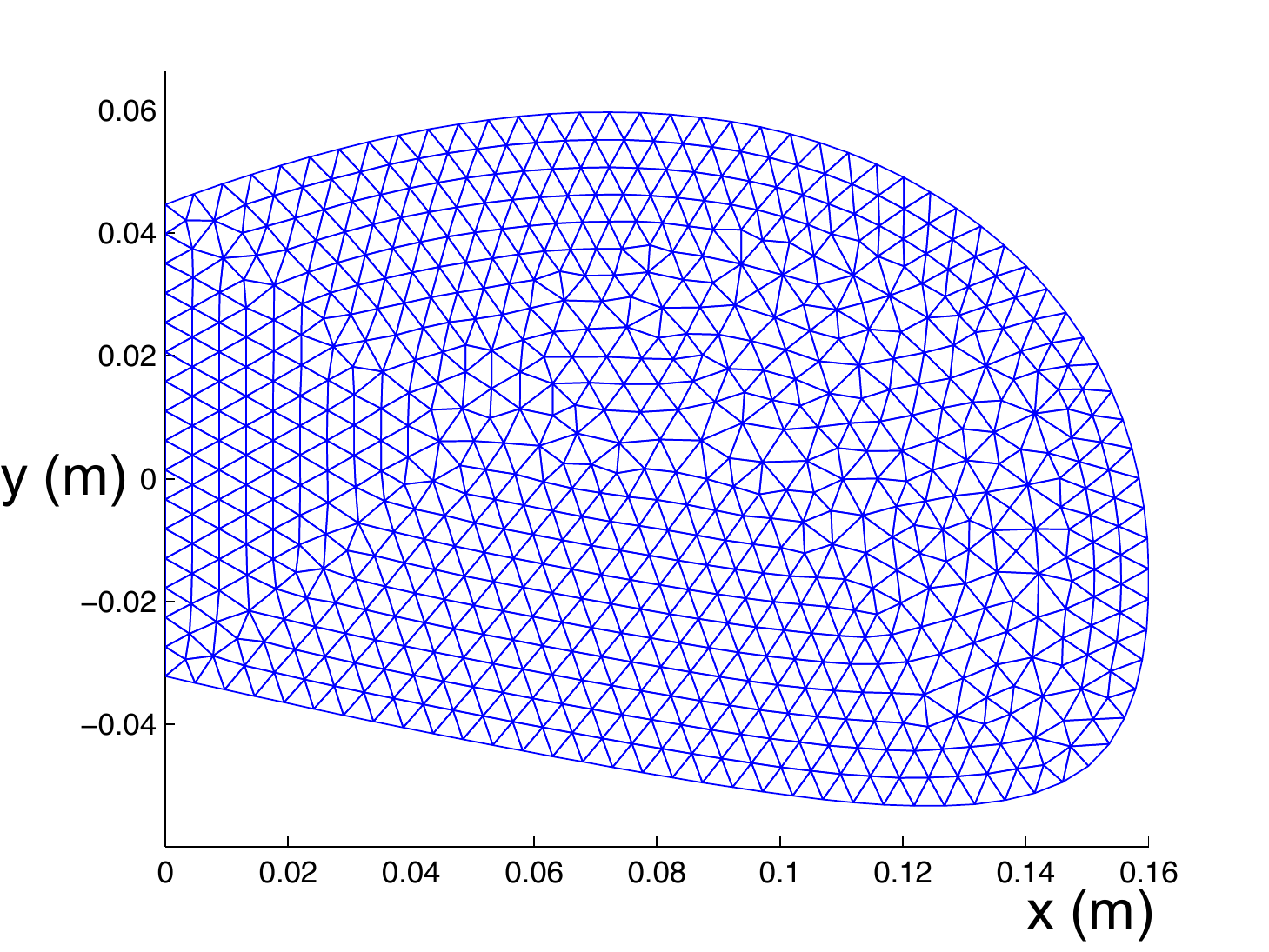}}
\par\end{centering}

\centering{}\subfloat[]{\begin{centering}
\includegraphics[scale=0.5]{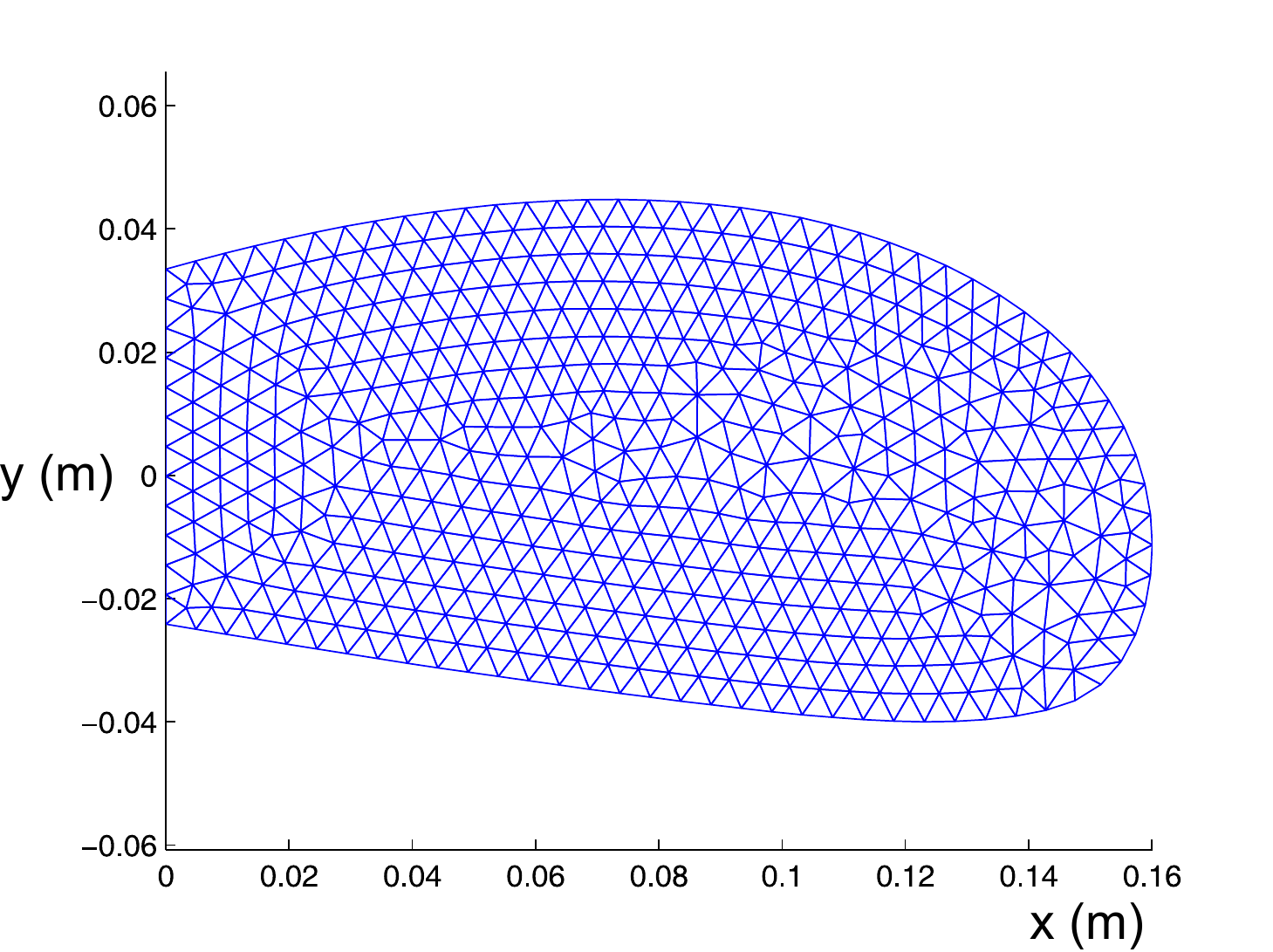}
\par\end{centering}

}\caption{\label{fig:T3-meshes-generated}T3 meshes generated for blade shapes;
(a) B5-75, (b) B5-60, (c) B5-45}
\end{figure}

The blades were chosen from the Wageningen-B five bladed series having
expanded area ratios (EAR) of 0.75, 0.6 and 0.45. The propellers were
chosen to have a diameter of 0.4 m with the hub (boss) having a diameter
of 0.08 m respecting the standards of Wageningen-B series. One special
characteristic of the Wageningen-B series is all propellers, apart
from 4-bladed propellers, have constant pitch distributions in the
radial direction \citep{Kuiper1992}, making the blades 2-dimensional
on the plane of the blade. Thus, 2-dimensional meshes were generated
for these three blade shapes using 3-noded triangular elements (Figure
\ref{fig:T3-meshes-generated}). The chosen pressure distributions
and required pitch to diameter ratios are given in Table \ref{tab:Pressure-variations}.
Further note that unlike a Controllable Pitch Propeller (CPP), a passive
shape adaptive propeller cannot change its pitch at the hub. Thus,
the pitch values presented in Table \ref{tab:Pressure-variations}
and throughout this paper are pitch values required at the tip of
the blade. 

\begin{table}
\begin{centering}
\begin{tabular}{ccccc}
\toprule 
P (kPa) & $\Delta$P (kPa) & P/D & $\phi$ (deg) & $\Delta\phi$ (deg)\tabularnewline
\midrule
\midrule 
180 & -70 & 0.7 & 12.56 & -3.44\tabularnewline
\midrule
\midrule 
205 & -45 & 0.8 & 14.3 & -1.7\tabularnewline
\midrule
\midrule 
250 (Cruise) & 0 & 0.9 & 16 & 0\tabularnewline
\midrule
\midrule 
270 & +20 & 1.0 & 17.66 & +1.66\tabularnewline
\midrule
\midrule 
300 & +50 & 1.1 & 19.3 & +3.3\tabularnewline
\bottomrule
\end{tabular}
\par\end{centering}

\caption{\label{tab:Pressure-variations}Pressure variations and required optimum
pitch values}
\end{table}

The optimization was carried out using prepreg AS4 carbon fiber reinforced
3501-6 epoxy with a nominal layer thickness of 125 \textgreek{m}m
(Table \ref{tab:Material-Properties} \citep{SodenProperties1998}).
It was considered that the blades were made of 40 layers. The lay
up was taken to be symmetric (i.e. - $\left[B\right]=\left[0\right]$)
to prevent warping during the manufacturing process. Thus, the optimization
algorithm was used to optimize 20 independent variables. It was initially
assumed that any fiber angle is possible to be manufactured; thus,
optimization was performed assuming ply angles could be any real number
(continuous variable ply optimization). Integer ply optimization was
also attempted and will be demonstrated in the proceeding section.
Ply angle results obtained using the optimization process, required
to enable the optimum pitch variation is given in Table \ref{tab:Ply-angle-results}
(all ply angle results specified in this paper are measured counter-clockwise
from x-axis to y-axis with coordinate system specified in Figure \ref{fig:T3-meshes-generated}). 

\begin{table}
\begin{centering}
\begin{tabular}{cc}
\toprule 
Property & Value\tabularnewline
\midrule
\midrule 
E\textsubscript{1}(GPa) & 126\tabularnewline
\midrule
\midrule 
E\textsubscript{2} (GPa) & 11\tabularnewline
\midrule
\midrule 
G\textsubscript{12} (GPa) & 6.6\tabularnewline
\midrule
\midrule 
\textgreek{n}\textsubscript{12} & 0.28\tabularnewline
\midrule
\midrule 
\textgreek{n}\textsubscript{23} & 0.4\tabularnewline
\midrule
\midrule 
Thickness (\textgreek{m}m) & 125\tabularnewline
\bottomrule
\end{tabular}
\par\end{centering}

\caption{\label{tab:Material-Properties}AS4-3501-6 Material Properties \citep{SodenProperties1998}}
\end{table}

\begin{center}
\begin{landscape}
\begin{table}
\begin{centering}
\begin{tabular}{cl}
\toprule 
Blade & \multirow{1}{*}{Ply angle configuration (deg) }\tabularnewline
\midrule
\midrule 
B5-75 & \multirow{1}{*}{{[}43.4/10.7/40.1/43.2/15.1/22.9/23.0/30.3/42.7/21.7/11.2/12.7/21.1/11.3/13.7/47.3/21.3/38.2/25.3/41.8{]}\textsubscript{s}}\tabularnewline
\midrule
\midrule 
B5-60 & {[}40.1/56.2/53.6/14.7/55.1/51.0/46.4/27.7/33.9/32.7/39.1/4.0/52.3/4.0/15.5/123.5/56.1/12.6/9.8/21.9{]}\textsubscript{s }\tabularnewline
\midrule
\midrule 
B5-45 & {[}11.7/53.6/28.9/51.6/37.8/28.2/36.6/51.0/51.4/31.0/24.7/13.9/21.8/3.6/35.2/52.0/26.5/26.9/44.2/39.4{]}\textsubscript{s} \tabularnewline
\bottomrule
\end{tabular}
\par\end{centering}

\caption{\label{tab:Ply-angle-results}Ply angle results}
\end{table}
 \end{landscape}
\par\end{center}

\begin{figure}
\begin{centering}
\includegraphics[scale=0.3]{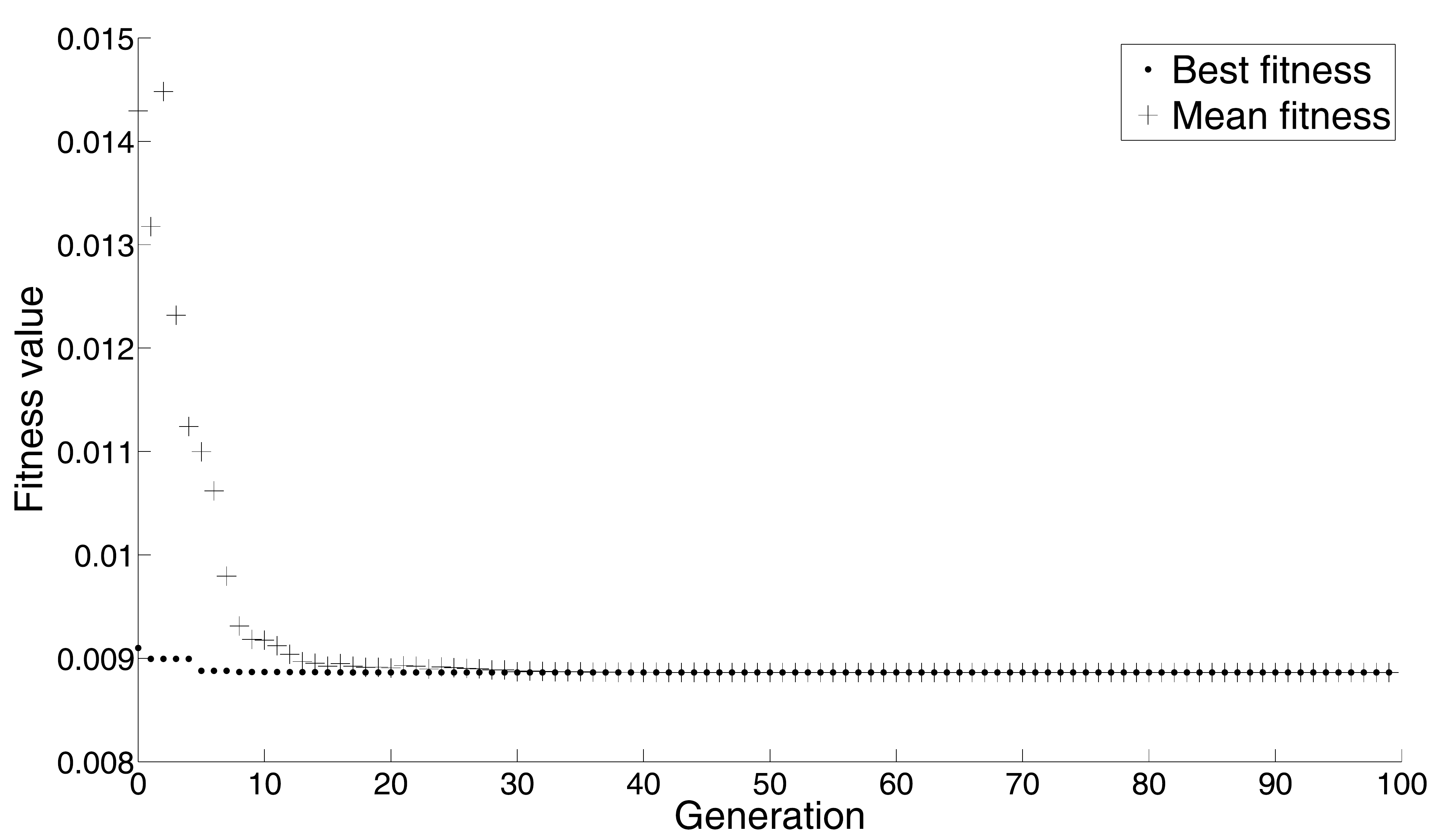}
\par\end{centering}

\caption{\label{fig:GA-Conv}Typical Convergence of the Genetic Algorithm}
\end{figure}

\begin{figure}
\begin{centering}
\includegraphics[scale=0.7]{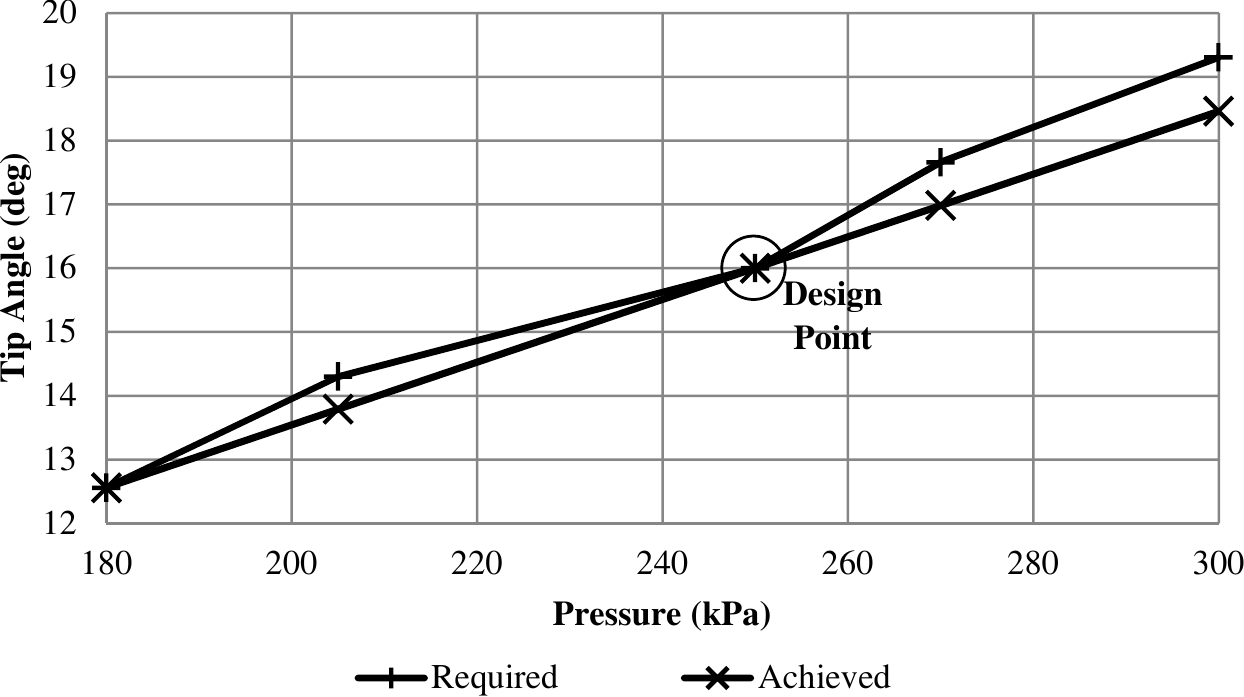}
\par\end{centering}

\caption{\label{fig:Comparison-between-pitch}Comparison between pitch angles}
\end{figure}

For the above three blade shapes the genetic algorithm converged to
a minimum objective function value of 0.51 deg ($8.86\times10^{-3}$rad).
Figure \ref{fig:GA-Conv} demonstrates the typical convergence of
GA for one of the optimization tasks. The minimum objective function
value can be thought of as the average deviation between the required
and the achieved result for each pitch angle. Thus, the results were
of high satisfaction. Although the requirement of tip angle vs pressure
variation was non-linear, the best a passive pitch adapting blade
can achieve is a linear tip angle variation with uniform pressure
as evident from Figure \ref{fig:Comparison-between-pitch}. 

Above ply angle results were verified using a commercial FE software
(ANSYS) and were found to provide bend-twist coupling characteristics
exactly as predicted by the in-house GA coupled FEM code. Thus, the
coupled code for GA and FEM can be used to optimize ply angles based
on the required deformation. Figure \ref{fig:Comparison-between-pitch}
demonstrates the comparison between achieved pitch angle and required
pitch angle. All three blade shapes behaved similarly and were able
to achieve the same pitch angles; thus, for brevity, only one comparison
is presented in Figure \ref{fig:Comparison-between-pitch}.

\subsection{Further analysis of Wageningen B5-45}

The Wageningen B5-45 blade was further analyzed for a weighted off-design
point case, integer optimization and the effect of layer thickness.
The unloaded shape required to behave according to Table \ref{tab:Pressure-variations}
was also obtained. For the weighted off-design point, it was assumed
that the $P=270\, kPa$ point has a relative weightage of 4, while
all other off-design points have a relative weightage of 1 (refer
to eq. \ref{eq:ObjFunc}). The value of 4 was chosen in this paper
purely for demonstration purposes; actual values of weightages have
to be evaluated for a particular propeller application based on its
probability to operate at each off-design point. Similar to previous
cases, continuous ply optimization was performed. Table \ref{tab:Weighted-Ply-Angle-results}
summarizes the obtained ply angle results and Figure \ref{fig:Wghted_Tip_Ang_Comp}
demonstrates and compares the pitch angle variation with equally weighted
optimization. Figure \ref{fig:Wghted_Tip_Ang_Comp} demonstrates that
assigning weightages enables the achieved pitch curve to move towards
more critical off-design point at the expense of deviating further
from other off-design points. This is further elaborated by the reduction
of difference between the achieved and required tip angle for $270\, kPa$
pressure (Figure \ref{fig:270_Error_Reduction}). Thus, weightages
have to be carefully chosen based on the operating conditions of the
propeller.

\begin{center}
\begin{landscape}
\begin{table}
\begin{centering}
\subfloat[]{\centering{}%
\begin{tabular}{|c|c|}
\hline 
B5-45 (Weighted) & {[}51.7/58.5/40.9/57.8/55.8/38.5/45.1/40.6/28.8/76.1/32.2/20.3/54.8/58.3/56.1/32.9/63.0/57.9/32.8/106.4{]}\textsubscript{s}\tabularnewline
\hline 
\end{tabular}}
\par\end{centering}

\centering{}\subfloat[]{\centering{}%
\begin{tabular}{|c|c|c|c|c|c|}
\hline 
Pressure (kPa) & 180 & 205 & 250 & 270 & 300\tabularnewline
\hline 
\hline 
Required $\phi$ (deg) & 12.56 & 14.3 & 16 & 17.66 & 19.3\tabularnewline
\hline 
Non-weighted $\phi$ (deg) & 12.56 & 13.79 & 16 & 16.98 & 18.46\tabularnewline
\hline 
Weighted $\left(w_{4}=4\right)$$\phi$ (deg) & 11.38 & 13.03 & 16 & 17.32 & 19.3\tabularnewline
\hline 
\end{tabular}}\caption{\label{tab:Weighted-Ply-Angle-results}Weighted ply angle results;
(a) Ply angles (degrees), (b) Resulting pitch angles}
\end{table}
\end{landscape}
\par\end{center}

\begin{center}
\begin{figure}
\begin{centering}
\includegraphics[scale=0.4]{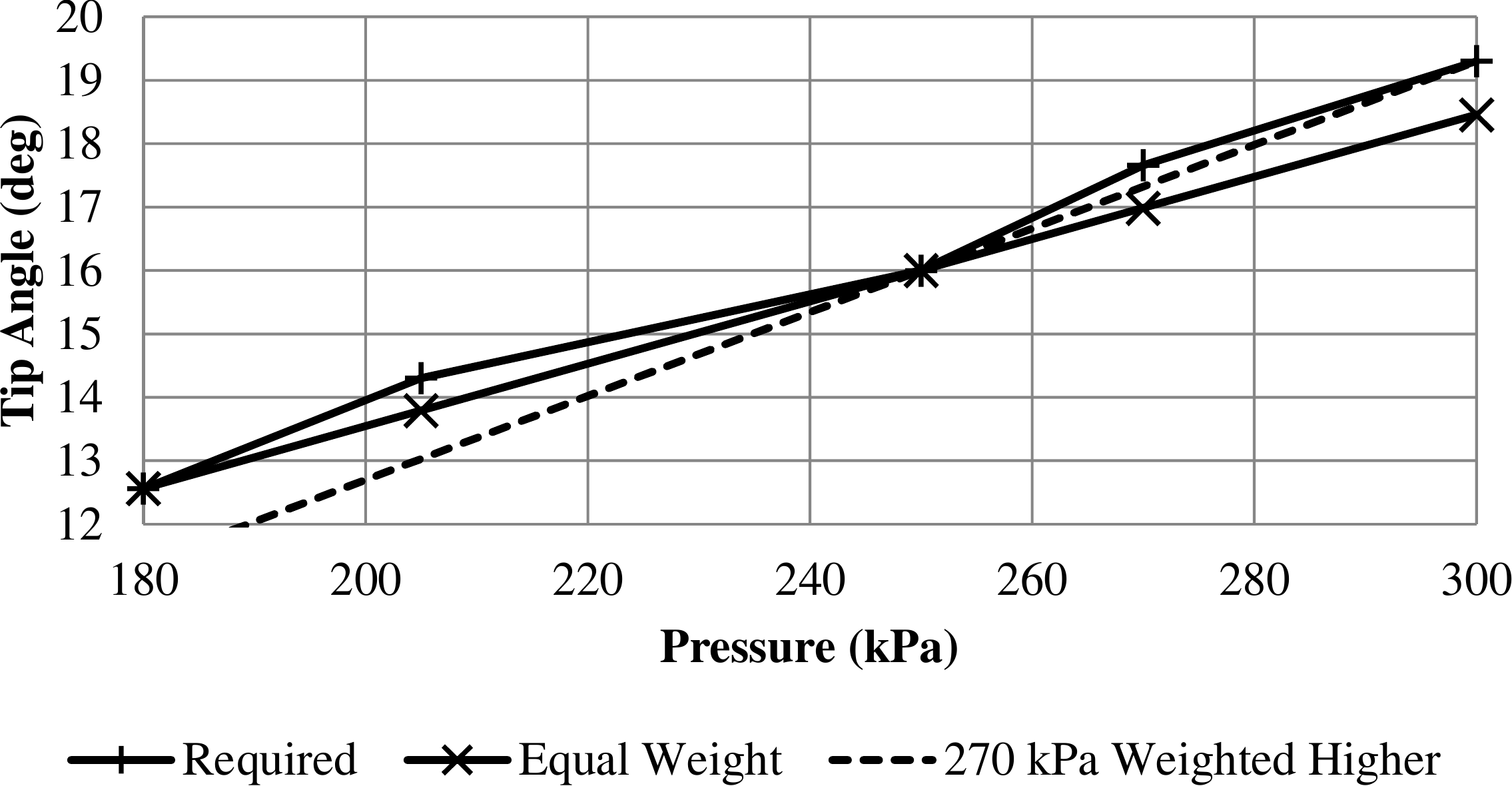}
\par\end{centering}

\caption{\label{fig:Wghted_Tip_Ang_Comp}Tip angle results for weighted optimization
case}
\end{figure}

\par\end{center}

\begin{center}
\begin{figure}
\begin{centering}
\includegraphics[scale=0.6]{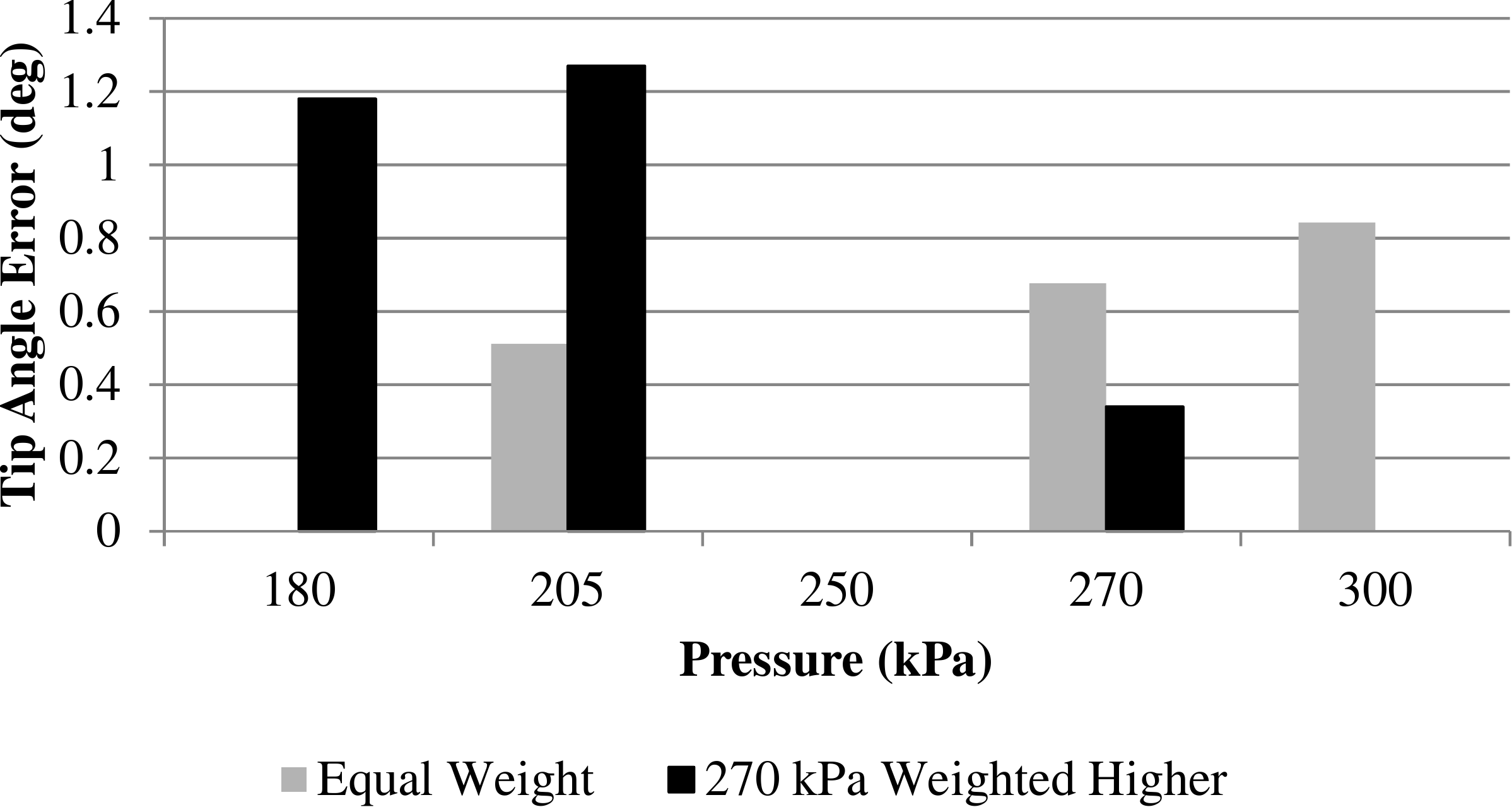}
\par\end{centering}

\caption{\label{fig:270_Error_Reduction}Difference between achieved and required
tip angles for non-weighted and weighted cases}

\end{figure}

\par\end{center}

Integer ply optimization was performed for equal weighting case. Integer
ply optimization is preferable as it greatly simplifies the lay up
and manufacturing process. It was observed that integer ply optimization
problems required a considerably longer time to converge at a result
compared to continuous ply optimization problems. Integer ply optimization
was performed for four cases: only integer ply angles $\left(0\leq\theta<\pi\right)$
are allowed, only products of 5 degrees are allowed, only products
of 10 degrees are allowed and only angles $\left\{ 0\text{\textdegree},30\text{\textdegree},45\text{\textdegree},60\text{\textdegree},90\text{\textdegree},120\text{\textdegree},135\text{\textdegree},150\text{\textdegree}\right\} $
are allowed. Results of these optimization problems are given in Table
\ref{tab:Integer-ply-results}. It was observed that for the first
three cases, the objective function was able to achieve the same minimum
value obtained for continuous ply optimization. Thus, limiting the
variable domain was not detrimental towards final results. The last
case with more aggressive variable domain limitations was found to
increase the value of the objective function. Hence, the achieved
tip angles further deviated from the required tip angles.

\begin{landscape}
\begin{table}
\begin{centering}
\begin{tabular}{ccc}
\toprule 
Case & Ply Angle Result (deg) & $f_{min}\:(rad)$\tabularnewline
\midrule
\midrule 
All Integer & {[}92/66/107/45/34/31/59/27/71/35/38/34/68/89/108/22/5/79/119/29{]}\textsubscript{s} & $8.86\times10^{-3}$ \tabularnewline
\midrule
\midrule 
Only 5n deg & {[}75/55/75/40/35/45/50/90/25/15/80/100/130/30/90/80/65/55/65/90{]}\textsubscript{s} & $8.86\times10^{-3}$ \tabularnewline
\midrule
\midrule 
Only 10n deg & {[}30/90/90/60/40/40/50/40/20/30/40/10/80/160/30/40/70/120/30/90{]}\textsubscript{s} & $8.86\times10^{-3}$ \tabularnewline
\midrule
\midrule 
Limited Ply Angles & {[}(60)\textsubscript{9}/90/60/90/150/60/135/(150)\textsubscript{5}{]}\textsubscript{s} & $4.41\times10^{-2}$ \tabularnewline
\bottomrule
\end{tabular}
\par\end{centering}

\caption{\label{tab:Integer-ply-results}Integer ply optimization results}
\end{table}
\end{landscape}

All above optimization tasks were carried out as unconstrained optimization.
One of the major drawbacks witnessed was the high flexibility of the
resulting blade, which resulted in a substantial rake deformation
in the process of pitch variation. It was investigated whether the
flexibility can be reduced by increasing the number of layers on the
blade. Inasmuch as, an attempt was made to optimize the blade assuming
there are 50 CFRP prepreg layers all having a nominal thickness of
$125\mu m$ (Table \ref{tab:Material-Properties}). It was observed
that, increasing the number of layers resulted in a stiffer blade,
which had relatively less bend-twist coupling performance compared
to the 40 layer blades attempted earlier (Table \ref{tab:different_thk_Results}).
However, the same twisting performance was able to be achieved if
the total laminate thickness of the 50 layer layup was the same as
the 40 layer layup; in other words, if the nominal layer thickness
of CFRP layers were reduced to $100\mu m.$ The resulting blade had
the same twisting performance as the 40 layer model, but with a less
lateral deformation. This was enabled by its higher bend-twist coupling
coefficient in the bending compliance matrix ($\left[d\right]$).
In addition, the layup with 40 $100\mu m$ layers had the same twisting
performance with a higher lateral bend. It must be noted that, layer
thickness is a practical limitation under the current composite prepreg
lamina manufacturing technology that is beyond composite designer's
control. However, it was clear that in order achieve best bend-twist
coupling performance one must attempt to reduce the the layer thickness
while increasing the number of layers in the laminate. Table \ref{tab:different_thk_Results}
summarizes results obtained for different layer numbers.

\begin{landscape}
\begin{table}
\begin{centering}
\begin{tabular}{c>{\raggedright}p{1.5cm}>{\centering}p{1.5cm}cc}
\toprule 
Layers & Layer Thk. ($\mu m$) & Total Thk. (mm) & Optimum Layup (deg) & $f_{min}\:(rad)$\tabularnewline
\midrule
\midrule 
\multirow{2}{*}{40} & \multirow{2}{1.5cm}{\centering{}125} & \multirow{2}{1.5cm}{\centering{}5} & {[}11.7/53.6/28.9/51.6/37.8/28.2/36.6/51.0/51.4/31.0/24.7/ & \multirow{2}{*}{0.00886}\tabularnewline
 &  &  & 13.9/21.8/3.6/35.2/52.0/26.5/26.9/44.2/39.4{]}\textsubscript{s} & \tabularnewline
\midrule
\midrule 
\multirow{2}{*}{50} & \multirow{2}{1.5cm}{\centering{}125} & \multirow{2}{1.5cm}{\centering{}6.25} & {[}41.9/41.9/78.8/41.6/43.5/42.3/45.1/79.4/42.8/41.6/78.8/44.8/44.4/ & \multirow{2}{*}{0.0134}\tabularnewline
 &  &  & 41.5/44.1/41.5/79.2/45.1/79.0/41.0/80.4/43.0/45.7/41.7/81.2{]}\textsubscript{s} & \tabularnewline
\midrule
\midrule 
\multirow{2}{*}{50} & \multirow{2}{1.5cm}{100} & \multirow{2}{1.5cm}{5} & {[}32.1/46.4/44.8/55.6/50.3/52.5/43.2/2.1/22.4/50.2/21.1/17.2/36.6/ & \multirow{2}{*}{0.00886}\tabularnewline
 &  &  & 36.5/13.5/40.8/38.1/49.3/45.3/26.9/52.3/12.2/55.8/8.8/3.3{]}\textsubscript{s} & \tabularnewline
\midrule
\midrule 
\multirow{2}{*}{40} & \multirow{2}{1.5cm}{100} & \multirow{2}{1.5cm}{4} & {[}4.1/30.3/11.3/7.6/24.3/109.4/60.4/0.41/42.9/15.7/15.2/ & \multirow{2}{*}{0.00886}\tabularnewline
 &  &  & 49.6/39.1/33.2/110.7/113.9/83.7/46.7/64.4/97.5{]}\textsubscript{s} & \tabularnewline
\midrule
\midrule 
\multirow{2}{*}{20} & \multirow{2}{1.5cm}{250} & \multirow{2}{1.5cm}{5} & \multirow{2}{*}{{[}15.7/19.7/30.3/49.7/37.2/44.3/22.1/34.9/47.0/3.4{]}\textsubscript{s}} & \multirow{2}{*}{0.00886}\tabularnewline
 &  &  &  & \tabularnewline
\bottomrule
\end{tabular}
\par\end{centering}

\caption{\label{tab:different_thk_Results}Results for different layer thicknesses}
\end{table}
 \end{landscape}

Note further that the current unconstrained optimization scheme is
capable of satisfying the pitch angle variation requirement almost
perfectly greatly enhancing the efficiency envelope. However, this
may not be possible due to the considerable increase in flexibility
and lateral deformations. Optimization constraints must be introduced
based on the application to limit such deformations at the expense
of efficiency gains.

\subsection{\label{sub:Unloaded-shape-calculation}Unloaded shape calculation}

The second stage of the proposed design process is the calculation
of the unloaded blade shape. The process is iteration based where
negative loads are applied to the shape at the design/cruise speed
of blade. Once the first unloaded shape is obtained, the loads at
design speed are applied to the unloaded shape to check whether it
reaches the desired shape at the operating speed. The desired shape
is the shape of the optimum alloy propeller chosen. If the shape requirement
has not been met, the initial shape is adjusted by the difference
and further iterated until there is satisfactory shape convergence
(Figure \ref{fig:Unloaded-shape-iteration}).

\begin{center}
\begin{figure}
\begin{centering}
\includegraphics[scale=0.7]{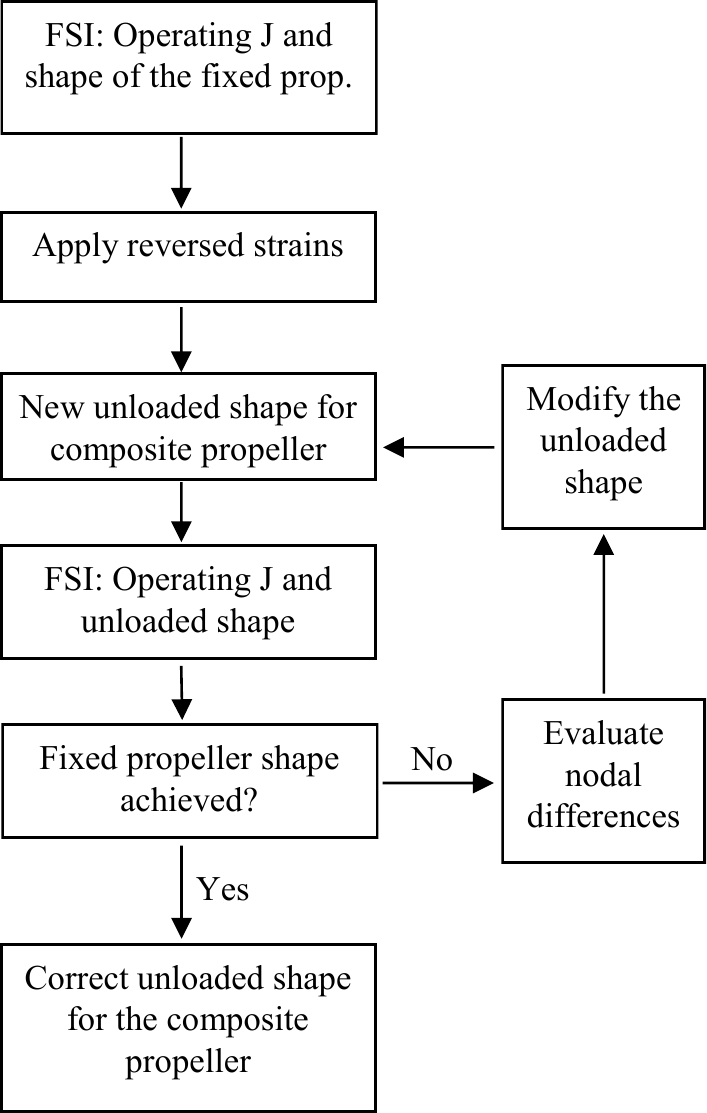}
\par\end{centering}

\caption{\label{fig:Unloaded-shape-iteration}Unloaded shape iteration process}
\end{figure}

\par\end{center}

\begin{figure}
\begin{centering}
\includegraphics[scale=0.4]{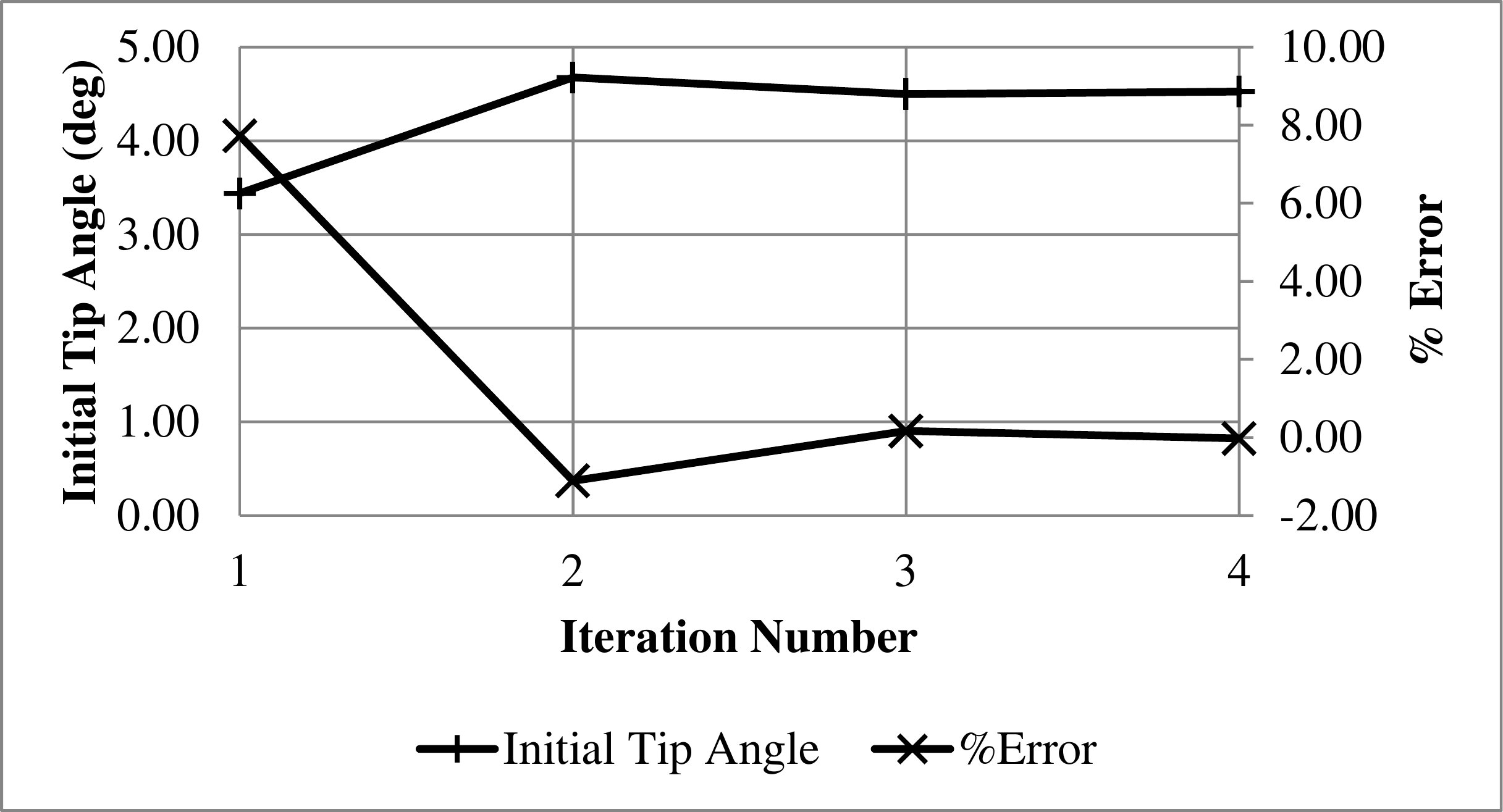}
\par\end{centering}

\caption{\label{fig:Unloaded-shape-Summary}Unloaded shape iteration summary}
\end{figure}

The unloaded shape for the B5-45 blade was calculated for the layup
given in Table \ref{tab:Ply-angle-results}. As summarized in Figure
\ref{fig:Unloaded-shape-iteration}, the first step of the process
was to apply reverse strains or reverse loads. In this example reverse
loads were applied. If reverse strains were applied the iteration
path would be different, but will eventually converge to the same
unloaded shape. Based on the iterations it was concluded that the
blade has to be manufactured with an initial uniform twist of 16\textsuperscript{o}
pitch angle at the root and 4.53\textsuperscript{o} at the tip. Further
loading due to the movement of the ship is expected to increase the
pitch of the propeller. Summary of iteration values is presented in
Table \ref{tab:Unloaded-shape-iteration}. These values are depicted
in Figure \ref{fig:Unloaded-shape-Summary} for better representation.

\begin{table}
\begin{centering}
\begin{tabular}{ccccc}
\toprule 
Iteration & Ini. Tip & Loaded Tip & \%Error & Adjustment\tabularnewline
\midrule
\midrule 
1 & 3.44\textsuperscript{o} & 14.76\textsuperscript{o}  & 7.73 & 1.24\textsuperscript{o}\tabularnewline
\midrule
\midrule 
2 & 4.68\textsuperscript{o}  & 16.18\textsuperscript{o}  & -1.11 & -0.18\textsuperscript{o}\tabularnewline
\midrule
\midrule 
3 & 4.50\textsuperscript{o}  & 15.97\textsuperscript{o}  & 0.17 & -0.03\textsuperscript{o}\tabularnewline
\midrule
\midrule 
4 & 4.53\textsuperscript{o} & 16.00\textsuperscript{o} & -0.03 & 0.00\textsuperscript{o}\tabularnewline
\bottomrule
\end{tabular}
\par\end{centering}

\caption{\label{tab:Unloaded-shape-iteration}Unloaded shape iteration steps}
\end{table}

\section{Conclusion}

The paper presented an optimization scheme for optimizing the ply
lay up of composite marine propeller blades using a coupled Genetic
Algorithm and Finite Element approach. The genetic algorithm was utilized
to vary the ply angles based on their fitness values while the Cell-Based
Smoothed FE code using Discrete Shear Gap Method, was used to evaluate
the deflection and pitch of the blade under given loadings. Pressure
distributions were chosen to be uniform and arbitrary, letting the
readers to use exact pressure distribution functions using CFD methods.
Continuous ply optimization was performed for three blade shapes to
demonstrates the design scheme applicability for various shapes. Further,
B5-45 was analyzed for weighted off-design point optimization, integer
ply optimization and effect of different layer thicknesses. Unloaded
shape for B5-45 was also calculated and was found to provide good
convergence. Results obtained using the coupled optimization algorithm
had a high accuracy and the process can be confidently used in future
research endeavors. The framework presented in the paper was purely
based on optimization for best efficiency using shape adaptability
achieved through optimization of composite layup angles. There are
many steps and adjustments to be made based on strength testing, cavitation
performance testing, natural frequency and structural divergence analysis,
etc. Thus, the proposed framework can be seen as a foundation for
designing shape-adaptive composite marine propellers.

\section*{Acknowledgments}

Authors would like to thank the Defence Science and Technology Organization
of Australia (DSTO) for the financial support and expertise provided
for the project. In addition, a special acknowledgment to Dr Andrew
Philips of DSTO for his constant support and input. S Natarajan would
like to acknowledge the financial support of the School of Civil and
Environmental Engineering, The University of New South Wales for his
research fellowship since September 2012.

\section*{References}

\bibliographystyle{elsarticle/elsarticle-num}
\bibliography{laminatedCompo}

\end{document}